\documentclass{dynafront2025} 

\usepackage[nosemicolon,noend]{algorithm2e}

\title[Automated algorithm design via NP interpolation]{Automated algorithm design via Nevanlinna-Pick interpolation }

\optauthor{%
\Name{Ibrahim K. Ozaslan$^1$} \Email{ozaslan@usc.edu}\\
\Name{Tryphon T. Georgiou$^2$} \Email{tryphon@uci.edu}\\
\Name{Mihailo R. Jovanovic$^1$} \Email{mihailo@usc.edu}\\
\addr $^1$Ming Hsieh Department of Electrical and Computer Engineering, University of Southern California, USA\\
$^2$Department of Mechanical and Aerospace Engineering, University of California, Irvine, USA}

\usepackage{float}
\usepackage{enumitem}

\makeatletter
\newcommand*\rel@kern[1]{\kern#1\dimexpr\macc@kerna}
\newcommand*\widebar[1]{%
	\begingroup
	\def\mathaccent##1##2{%
		\rel@kern{0.8}%
		\overline{\rel@kern{-0.8}\macc@nucleus\rel@kern{0.2}}%
		\rel@kern{-0.2}%
	}%
	\macc@depth\@ne
	\let\math@bgroup\@empty \let\math@egroup\macc@set@skewchar
	\mathsurround\z@ \frozen@everymath{\mathgroup\macc@group\relax}%
	\macc@set@skewchar\relax
	\let\mathaccentV\macc@nested@a
	\macc@nested@a\relax111{#1}%
	\endgroup
}

\DeclareMathOperator*{\minimize}{minimize}

\DeclareMathOperator*{\subjectto}{subject~to}

\newcommand{\reals}[1]{\mathbb{R}^{#1}}
\newcommand{\cN}{\mathcal{N}}
\DeclareMathOperator{\re}{Re}

\DeclareMathOperator{\diag}{\mathrm{diag}}

\newcommand{\norm}[1]{\| #1 \|}    
\newcommand{\inner}[2]{\langle #1,#2\rangle}    

\newcommand{\ds}{\displaystyle}
\newcommand{\myset}[1]{\{ #1 \}}
\newcommand{\enma}[1]{\ensuremath{#1}}

\newcommand{\non}{\nonumber}

\newcommand{\beq}{\begin{equation}}
\newcommand{\eeq}{\end{equation}}
\newcommand{\ba}{\begin{array}}
\newcommand{\ea}{\end{array}}
\newcommand{\bseq}{\begin{subequations}}
\newcommand{\eseq}{\end{subequations}}

\newcommand{\DefinedAs}[0]{\mathrel{\mathop:}=}

\newcommand{\rme}{\mathrm{e}}

\newcommand{\matbegin}{\left[}
\newcommand{\matend}{\right]}

\newcommand{\tbo}[2]{
		\matbegin \begin{array}{c}
				#1 \\ #2
			\end{array} \matend }

\newcommand{\obt}[2]{
	\matbegin \begin{array}{cc}
		#1 & #2
	\end{array} \matend }

\newcommand{\tbt}[4]{
		\matbegin \begin{array}{cc}
				#1 & #2 \\ #3 & #4
			\end{array} \matend }

\newcommand{\ww}{w}

\newcommand{\sig}{\sigma}

\newcommand{\xs}{x^\star}
\newcommand{\ws}{\ww^\star}

\newcommand{\dds}{\delta^\star}

\newcommand{\xb}{\widebar{x}}

\newcommand{\hb}{\widebar{h}}

\newcommand{\xh}{\widehat{x}}

\newcommand{\gradh}{\widehat{\nabla f}}

\newcommand{\vh}{\widehat{v}}

\newcommand{\cZ}{\enma{\mathcal Z}}

\newcommand{\cL}{\enma{\mathcal L}}

\newcommand{\cK}{\enma{\mathcal K}}

\newcommand{\cB}{\enma{\mathcal B}}
\newcommand{\cH}{\enma{\mathcal H}}

\newcommand{\cA}{\enma{\mathcal A}}

\newcommand{\sigu}{\widebar{\sigma}}
\newcommand{\sigl}{\underline{\sigma}}

\newcommand{\cHb}{\enma\widebar{\mathcal H}}
\newcommand{\Deltab}{\enma\widebar{\Delta}}

\newcommand{\bD}{\enma{\mathbb D}}

\newcommand{\Rhardy}{\enma{\mathcal{RH}}_\infty}

\begin{document}

\maketitle

\begin{abstract}%
The synthesis of optimization algorithms typically
follows a design-first-analyze-later approach, which often obscures fundamental performance limitations and hinders the
systematic design of algorithms operating at the achievable theoretical boundaries. Recently, a framework based on frequency-domain techniques from robust control theory has emerged
as a powerful tool for automating algorithm synthesis. By
integrating the design and analysis stages, this framework
enables the identification of fundamental performance limits.
In this paper, we build on this framework and extend it to
address algorithms for solving strongly convex problems with equality constraints. As a result, we obtain a new class of algorithms that offers sharp trade-off between number of matrix multiplication per iteration and convergence rate.
\end{abstract}

\vspace*{-1ex}
\section{Introduction}
	\label{sec.intro}

The design and analysis of optimization algorithms typically begins with selecting an appropriate algorithmic structure based on the problem's optimality conditions. This is followed by a convergence analysis, which often relies on creative and problem-specific reasoning. Such a {\em design-first-analyze-later\/} approach typically requires deriving tight analytical inequalities$-$an inherently challenging and non-systematic task. This difficulty contributes, in part, to the limited understanding of fundamental performance limits within specific algorithm classes.

Viewing optimization algorithms through the lens of dynamical systems and applying Lyapunov stability theory enables a system-theoretic understanding of their behavior. The challenging aspect comes from the fact that designing Lyapunov functions that lead to the optimal convergence rates requires deep problem-specific insight~\cite{suboycan16,wibwiljor16,hules17,shidujorsu21,muejor19, ozajovAUT25,zhezhusoblali23}. A complementary approach involves fixing the structure of the Lyapunov function and the algorithm, and then using computational tools to search for valid certificates and parameters~\cite{lesrecpac16, vanfrelyn17, micschebe21,schebehol23,boyparryusuh24}. While this method has the potential to yield sharp results~\cite{vanfrelyn17, schebe21}, its reliance on fixed functional forms restricts its capacity to uncover fundamental performance limits.

Recently a frequency-domain approach leveraging tools from robust control theory emerged as a powerful tool for automating algorithm design and analysis~\cite{ugrpetsha23,zhawulichegeo24,wuchejovgeo24}. This framework models optimization algorithms as transfer functions within a feedback system, derives analytic conditions to guarantee desirable algorithmic properties, and uses Nevanlinna-Pick interpolation to construct transfer functions that meet these conditions. When the algorithm specifications are expressed through necessary and sufficient conditions, this framework provides a means for identifying fundamental performance limits. For instance, a constructive proof that establishes Nesterov's lower bound on the complexity of gradient-based methods in unconstrained minimization, which was originally obtained for infinite dimensional problems using adversarial examples, has been recently provided in~\cite{ugrpetsha23} using this framework, and then independently in~\cite{wuchejovgeo24} building on~\cite{zhawulichegeo24}.

In this paper, we extend this framework to address algorithms for convex optimization problems with equality constraints,
	\begin{align}
	\label{eq.main}
	\minimize\limits_{x} ~f(x)
~
	\subjectto ~Ex \,=\, q, ~~ E\in\reals{d\times n}
	\end{align}
where $f$ is an $m$-strongly convex function with an $L$-Lipschitz continuous gradient. Such problems arise in a host of applications spanning federated learning~\cite{lisahtalsmi20}, imaging~\cite{chapoc16}, optimal transport~\cite{peycut19}, and fluid dynamics~\cite{tahgonsho23}. A particularly important instance is the decentralized optimization where $E$ is a gossip (or Laplacian) matrix of the network, $q=0$, and the objective function is the sum of local costs $f(x) = \sum_{i=1}^Nf_i(x_i)$~\cite{scabacbubleemas17}.

\paragraph{Existing algorithms:} A lower bound on the iteration complexity of the gradient-based algorithms for problem~\eqref{eq.main}, along with an algorithm achieving the bound, was established in~\cite{salconkovric22}. This algorithm attains $\epsilon$-accuracy in $O(\sqrt{\kappa_f}\log\tfrac{1}{\epsilon}))$ iterations, but requires $O(\sqrt{\kappa_E})$ matrix-vector multiplications per iteration. Here, $\kappa_f \DefinedAs L/m$ and $\kappa_E \DefinedAs \sig_1/\sig_r$, where $\sig_1$ and $\sig_r$ denote the largest and smallest (nonzero) singular values of $E^\top E$, respectively. When matrix-vector multiplications dominate the cost, e.g., in distributed optimization with limited communication budget~\cite{lanleezho20, quli19}, inner-loop iterations are undesirable. Among single-loop methods, Gradient Descent-Ascent (GDA) achieves $O(\kappa_f\kappa_E\log\tfrac{1}{\epsilon})$ iteration complexity, while a variant of Condat-Vu algorithm~\cite{con13,vu13}, known as the PAPC~\cite{salconmisric22}, achieves $O(\max(\kappa_f,\kappa_E)\log\tfrac{1}{\epsilon})$. In the special case, where $E$ is the gossip matrix and $q=0$, the accelerated PAPC algorithms attains $O(\max(\sqrt{\kappa_f\kappa_E},\kappa_E)\log\tfrac{1}{\epsilon}))$ complexity~\cite{kovsalric20}, while accelerated EXTRA achieves $O(\sqrt{\kappa_f\kappa_E}\log(\kappa_f\kappa_E)\log\tfrac{1}{\epsilon})$~\cite{lilin20}.

\vspace*{1ex}
\noindent\textbf{Our contribution} lies in the general framework we develop, which potentially enables automated algorithm design across diverse settings. Using this framework, we synthesize the Interpolated Gradient Method (I-GM) for solving problem~\eqref{eq.main} with an elegant structure that provides a sharp trade-off between the number of matrix-vector multiplications per iteration, $\ell$, and the worst-case convergence rate. The resulting iteration complexity is $O(\max(\kappa_f,\kappa_E/\ell)\log\tfrac{1}{\epsilon})$. For example, in the decentralized setting with limited computation/communication budget, one can effectively balance the local information propagation and convergence rate by adjusting the parameter $\ell$.

\vspace*{-1ex}
\section{Automated algorithm synthesis}\label{sec.framework}  
Our framework consists of the following steps.  We first outline the algorithm specifications (Sec.~\ref{sec.specs}) and utilize $\cZ$-transform to represent the class of algorithms under consideration (Sec.~\ref{sec.representation}). Then, we formulate the algorithm design as a controller synthesis problem (Sec.~\ref{sec.controller}) and obtain analytical characterization of the design specifications (Sec.~\ref{sec.characterization}) in terms of the underlying transfer functions. By recasting control synthesis as an interpolation problem (Sec.~\ref{sec.interpolation}), we provide a solution that leads to an implementable optimization algorithm with a guaranteed convergence rate (Sec.~\ref{sec.algo}). 

For clarity of exposition, we simplify our synthesis procedure by replacing  the constraint with $E^\top Ex = E^\top q$ and assuming $q=0$. The first modification incurs no loss of generality, as multiplying the equality constraint by $E^\top$ leaves the feasible set unchanged. Crucially, the synthesized algorithm never explicitly forms $W \DefinedAs E^\top E$, it requires only matrix-vector products with $E^\top$ and $E$. The second assumption, $q=0$, is purely expository; extending our framework to handle arbitrary $q$ is achieved by a straightforward modification in Section~\ref{sec.algo}  while preserving all theoretical guarantees of the synthesized algorithm.

\subsection{Design specifications}\label{sec.specs}
\paragraph{Linearity:} 
We investigate optimization algorithms with linear updates both in optimization variables and gradients.
\vspace*{-1ex}
\paragraph{Explicitness:} 
We design single-loop algorithms that avoid circular dependence in the evaluation of $x^{k+1}$. In other words, the computation of $x^{k+1}$ should not exploit $Ex^{k+1}$ or $\nabla f(x^{k+1})$. This requirement can be relaxed if one can afford to utilize linear system solvers~\cite{ozapilari23} or evaluation of proximal operators~\cite{wuchejovgeo24} in the iterative scheme.

\vspace*{-1ex}
\paragraph{Optimality:}
The algorithm should asymptotically converge to the unique solution $\xs$ of problem~\eqref{eq.main}. A necessary and sufficient condition for $\xs$ (when $q=0$) is given by the existence of unique vectors $\dds_1\in\reals{r}$ and $\ws_2\in\reals{n-r}$ such that
$\nabla f(\xs)  = V_1\dds_1$ and $\xs  =  V_2\ws_2$. Here, $r<n$ determines the rank of matrix $W = E^\top E$ with the singular value decomposition,
\beq\non
W \,=\, \obt{V_1}{V_2}\tbt{\Sigma}{0}{0}{0}\tbo{V_1^\top  }{V_2^\top  },~~ \Sigma \,\DefinedAs\, \diag(\sig_1,\dots,\sig_r), ~~ V \,\DefinedAs\, \obt{V_1}{V_2}.
\eeq

\vspace*{-1ex}
\paragraph{Linear convergence:}
We require algorithms to linearly converge at a rate of at least $\rho$, i.e., there exist $\rho\in(0,1)$ and $M>0$ such that $\norm{x^k - \xs}   \leq  M\norm{x^0 - \xs}\rho^k$ for any $x^0\in\reals{n}$.

\subsection{Algorithm representation}\label{sec.representation}

 The algorithms with linear updates both in optimization variables and gradients take the form of
\beq\label{eq.repre}
z\cK_0(z,W)\xh(z)  \,=\, \cK_1(z,W)\xh(z)   \,+\, \cK_2(z,W)\gradh(z)
\eeq
where $\xh(z)$ and $\gradh(z)$ are the $\cZ$-transforms of sequences $\{ x^k \}_{k=0}^\infty$ and $\{ \nabla f(x^k) \}_{k=0}^\infty$. Here, transfer function $\cK_0(z,W)$ determines the mapping between the optimization variable and the input of the gradient, while $\cK_1(z,W)$ and $\cK_2(z,W)$ determine how the past iterates and past gradients, respectively, are used to update the optimization variable. These three transfer functions are (potentially) parametrized by constraint matrix $E$, or in this particular setting, by $W \DefinedAs E^\top E$.

The class of algorithms that can be brought into the form of~\eqref{eq.repre} includes, but not limited to, accelerated methods such as Nesterov's Accelerated Gradient~\cite{nes03}, distributed algorithms such as DIGing~\cite{nedolsshi17}, EXTRA~\cite{shilinwuyin15}, and $\cA\cB\cN$~\cite{xinjakkha19} as well as descent-ascent type methods such as GDA~\cite{algsay20} and PAPC method~\cite{salconkovric22}; see appendix for a detailed example.

\subsection{From algorithm to controller design}\label{sec.controller}

\begin{figure*}[h]
	\centering
	\includegraphics[width=1\textwidth]{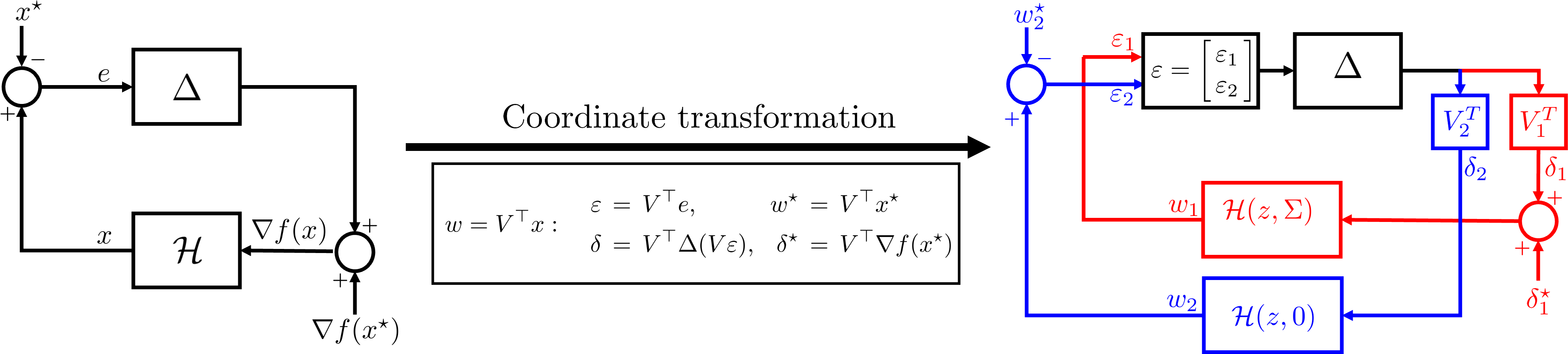}
	\caption{(Left) Lur'e system representing the optimization algorithms. (Right) Decomposition of the system on the left into two subsystems coupled by the nonlinear element.}
	\label{fig.diagram}
\end{figure*}

The optimization algorithms given by~\eqref{eq.repre} can be viewed as the Lur'e system in Fig.~\ref{fig.diagram}, that is, a feedback interconnection between an LTI system defined by input-output map $\xh(z) = \cH(z,W)\gradh(z)$ and a static nonlinear block $\Delta$,  where
\beq\label{eq.lure}
\cH(z,W)  \,\DefinedAs\, (z\cK_0(z,W) \,-\, \cK_1(z,W))^{-1}\cK_2(z,W), ~~ \Delta( e )  \,\DefinedAs\, \nabla f(e  + \xs)  \,-\, \nabla f(\xs).
\eeq
The nonlinear map $\Delta$ is sector bounded in $[m,L]$, i.e.,
$
m\norm{e}^2 \leq \inner{\Delta(e)}{e}  \leq L\norm{e}^2$ for all  $e \in\reals{n}
$.

Our task is to design a linear controller $\cH(z,W)$, representing the optimization algorithm, that satisfies design specifications, optimality, explicitness, and linear convergence (the linearity is enforced by formulation~\eqref{eq.repre}). To this end, we obtain an analytical characterization of design specifications, but first, we employ the following transformations to bring the Lur'e system into a form that is more convenient for analysis/synthesis. 


\vspace*{-1ex}
\paragraph{Coordinate transformation:}

We apply the change of  coordinates $w = V^\top x$ using the singular vectors of $W$ (see Fig.~\ref{fig.diagram} for additional definitions) to  cast transfer function $\cH$ as, 
\beq\non
\cH(z,W)  \,=\,  \obt{V_1}{V_2}\tbt{\cH(z, \Sigma)}{0}{0}{\cH(z, 0)}\tbo{V_1^\top  }{V_2^\top  }.
\eeq 
Since the resulting subsystems are in the diagonal form, they can be further split into the scalar systems, each of which shares the same scalar transfer function parameterized by either singular values or zero, i.e., $\cH(z,\Sigma)   =  \diag\big(\{h(z,\sig_i)\}_{i=1}^r\big)$ and $\cH(z,0)   \,=\,  h(z,0)I_{n-r}$.

\vspace*{-1ex}
\paragraph{Loop transformation:}
We apply loop transformation~\cite[Sec.\ 6.5]{kha02} to the Lur'e system, enabling the use of standard passivity-based stability criteria without introducing additional conservatism. The transformation maps the pair $(\cH, \Delta)$ with sector-bounded $\Delta$ in $[m,L]$ to $(\cHb, \Deltab)$ with transformed nonlinearity $\Deltab$ being sector-bounded in $[0,\infty)$.  Using the coordinate transformation, we can also diagonalize the transformed linear system with the following one-to-one relation,
\beq\label{eq.onetoone}
h(z,\sigma)  \,=\, (\hb(z,\sigma) \,-\, 1)\big/(m\hb(z,\sigma) \,-\, L).
\eeq

Our task is now reduced to designing a scalar transfer function $\hb(z,\sigma)$ that satisfies the design specifications for all $\sigma \in [\sigl,\sigu]\cup\myset{0}$ where $\sigu\geq\sig_1\geq\sig_r\geq\sigl>0$.

\subsection{Analytical characterization of design specifications}\label{sec.characterization}

\paragraph{Explicitness as causality:}
A necessary condition for explicitness is the strict causality of the linear controller in the Lur'e system. The controller is strictly causal when transfer function $\cH(z,W)$ is strictly proper, meaning $\cH(\infty,W) = 0$. This is equivalent to $\cHb(\infty, W)  =  I$, i.e.,
	\bseq
	\label{eq.interpolation}
	\beq\label{cond.causal}
	\hb(\infty,\sigma)  \,=\,  1,
	~~\forall \, \sigma \,\in\,[\sigl,\sigu]\,\cup\,\myset{0}. 
	\eeq

\paragraph{Optimality as input tracking and disturbance rejection:}
The asymptotic convergence of the optimization algorithm to the optimal solution is equivalent, in the Lur'e system, to the asymptotic tracking of the step input $\tfrac{z}{z-1} \xs$ with unknown magnitude $\xs$ and rejection of the step disturbance $\tfrac{z}{z-1} \nabla f(\xs)$ with unknown magnitude $\nabla f(\xs)$. The internal model principle~\cite{frawon76} implies that asymptotic tracking and rejection can be achieved if the subsystems $\cH(z,\Sigma)$ and $\cH(z,0)$ have a blocking zero and a pole at $z=1$, respectively. For the transformed system, this is equivalent to
\beq\label{cond.opt}
\hb(1,\sig) = 1,~~  \forall \sig \,\in\, [\sigl,\sigu], ~~~~~ \hb(1,0) \,=\, L/m .
\eeq

\paragraph{Linear convergence as $\rho$-stability:}
A necessary and sufficient condition for linear convergence of the algorithm is the $\rho$-stability of the Lur'e system. The circle criterion~\cite[Sec.\ 7.1]{kha02} implies $\rho$-stability of the transformed Lur'e system $(\cHb, \Deltab)$, if the transfer function $\hb(\gamma z,\sigma)$ is {\em strictly positive real\/} for all $\gamma\in(\rho, 1]$. This can be equivalently expressed as, for all $\gamma\in(\rho, 1]$ and $\sigma \in[\sigl,\sigu]\cup\myset{0}$,
\beq\label{cond.stab}
\ba{rl}
\mathrm{i})~~\hb(\gamma z,\sigma)\in\Rhardy,
~~~~~
\mathrm{ii})~~\re \, ( \hb(\gamma z,\sig) ) > 0, ~\forall z \in \partial\bD
\ea
\eeq
\eseq
where $\Rhardy$ denotes set of all real, rational, stable transfer functions,  $\re (\cdot)$ is the real part of a complex number, $\partial\bD\DefinedAs\myset{z\in\mathbb{C}, \, |z|=1}$, and $\mathbb{C}$ is the complex space. 


\vspace*{-1ex}
\subsection{From controller design to interpolation} \label{sec.interpolation}
With the algorithm representation~\eqref{eq.repre} and analytic characterization of design specifications~\eqref{eq.interpolation}, the controller design task reduces to the following algebraic problem:
\beq
\tag{IP}
\label{eq.IP}
\text{Find a $\rho\in(0,1)$ and design a $\hb(z,\sig)$ such that conditions~\eqref{eq.interpolation} hold for every $\gamma\in(\rho,1]$.}
\eeq
Using suitable conformal maps, we can cast~\eqref{eq.IP} as the following interpolation problem:
\beq
\non
\label{eq.interpolation2}
\ba{c}
\text{Given data pairs $\myset{(z_i, w_i)}_{i=1}^N$ with each $z_i\in\bD$, find a holomorphic function $\psi(z)$}
\\
\text{that satisfies $\{ |\psi(z)|\leq 1,$ for all $z\in\bD \}$ and $\psi(z_i) = w_i$, where $ \bD\DefinedAs \myset{z\in\mathbb{C},\,|z|< 1}$.}
\ea
\eeq 
The interpolation problem is then solved using the Nevanlinna-Pick (NP) method~\cite[Ch.\ 9]{doyfratan13}, which determines: (i) when a solution exists, i.e., a lower bound on $\rho$, via the Pick matrix; and (ii) how to parametrize all solutions, i.e., all algorithms satisfying design specifications, through an iterative construction based on Schur's algorithm. 

For any given constants $\ell>0$, $\alpha_1\in(0, 2/(L+m))$, $\alpha_2\in(0, 1/\sigu)$, and $g(\sig) \DefinedAs (1 - \alpha_2\sig)^\ell$, the NP method yields the following solution to problem~\eqref{eq.IP},
\beq \label{eq.theorem}
\rho \,=\, 
\max\big(1-\alpha_1m, \, (1-\alpha_2\sigl)^\ell\big), 
~~~~
\hb(\gamma z,\sigma) \,=\, \dfrac{\gamma z(\gamma z-g^2(\sig)) \,+\, \rho g^2(\sig)(\gamma z-1)}{\gamma z(\gamma z-g^2(\sig)) \,-\, \rho g^2(\sig)(\gamma z-1)}.
\eeq
The convergence rate $\rho$ is optimized by setting $\alpha_1 = 2/(m+L)$ and $\alpha_2=1/\sig_1$, resulting in $\rho^\star = \max\big(1 - 2/(\kappa_f+1), \boldsymbol{(}1-1/\kappa_E\boldsymbol{)}^\ell\big)$.

\vspace*{-1ex}
\subsection{Algorithm synthesis: Interpolated Gradient Method (I-GM)}\label{sec.algo}
The transfer functions in~\eqref{eq.repre} follow from $\cH$ in~\eqref{eq.lure} which itself is obtained via the one-to-one relation~\eqref{eq.onetoone} from~\eqref{eq.theorem}, leading to Algorithm~\ref{algo.1}.

When $q\neq 0$, the solution to problem~\eqref{eq.main} decomposes as $\xs = \xs_{\shortparallel} + \xs_\perp$ where $x_{\shortparallel}$ is the particular solution to $E^\top E x = E^\top q$ and $x_\perp$ is the null space component. Hence, to preserve the orthogonality between the two exogenous inputs, we replace $\xs$ in Fig.~\ref{fig.diagram} with $\xs_\perp$ and apply the substitution $E^\top Ex\leftarrow E^\top E(x - \xs_{\shortparallel} )$ in $\cH$. This ensures that the output of the LTI system is $x - \xs_{\shortparallel} $ and the error takes the correct form: $e= x -  \xs_{\shortparallel} - \xs_\perp = x - \xs$. In summary, substitution of  $E^\top Ex\leftarrow E^\top (Ex - q)$  in Algorithm~\ref{algo.1} yields Algorithm~\ref{algo.2} with $p_\ell(\sig) \DefinedAs \sum_{i=0}^{2\ell-1}(1 - \alpha_2\sig)^i$ (or $g(\sigma)= 1 - \alpha_2\sig p_\ell(\sig)$). 

\vspace{1ex}
\begin{minipage}{0.485\textwidth}
	\hspace*{-4ex}
	\SetAlgorithmName{Algorithm}{}{}
	\begin{algorithm2e}[H]
	\caption{ I-GM ($q=0$)}
	\SetAlgoLined
	\KwData{\small $v^{-1}\! = \!x^0$, $\alpha_1 \!\in\! (0, \tfrac{1}{L})$, $\alpha_2 \!\in \!(0, \tfrac{1}{\sigu})$, $\ell\!>\! 0$ }
	\For{$k = 0, 1, 2, \ldots$}{
		$v^k = x^k - \alpha_1\nabla f(x^k)$\;
		$x^{k+1} = (I - \alpha_2E^\top E)^{2\ell}(x^{k} + v^{k} - v^{k-1})$\;
	}
	\label{algo.1}
	\end{algorithm2e}
\end{minipage}\hspace*{-1ex}%
\begin{minipage}{0.5\textwidth}
	\centering
		\SetAlgorithmName{Algorithm}{}{}
	\begin{algorithm2e}[H]
	\caption{I-GM}
	\SetAlgoLined
	\KwData{\small $v^{-1}\! = \!x^0$, $\alpha_1 \!\in\! (0, \tfrac{1}{L})$, $\alpha_2 \!\in \!(0, \tfrac{1}{\sigu})$, $\ell\!>\! 0$ }
	\For{$k = 0, 1, 2, \ldots$}{
		$v^k = x^k - \alpha_1\nabla f(x^k)$\;
		$w^{k} = x^{k} + v^{k} - v^{k-1}$\;
		$x^{k+1}\! =\! w^{k} \!- \alpha_2\, p_\ell(E^\top E)(E^\top Ew^{k} \!-\! E^\top q)$\;
	}
	\label{algo.2}
	\end{algorithm2e}
\end{minipage}

\appendix









\newpage
\clearpage
\appendix
\section{Computational Experiments}

\begin{figure*}[h]
		\centering
		\begin{tabular}{c@{\hspace{-0.3 cm}}c@{\hspace{0.1 cm}}c@{\hspace{-0.3 cm}}c}
			\begin{tabular}{c}
				\vspace{.25cm}
				\normalsize{\rotatebox{90}{ $\log_{10} \norm{x^k - \xs}/\norm{\xs}$}}
			\end{tabular}
			&
			\begin{tabular}{c}
				\includegraphics[width=0.4\textwidth]{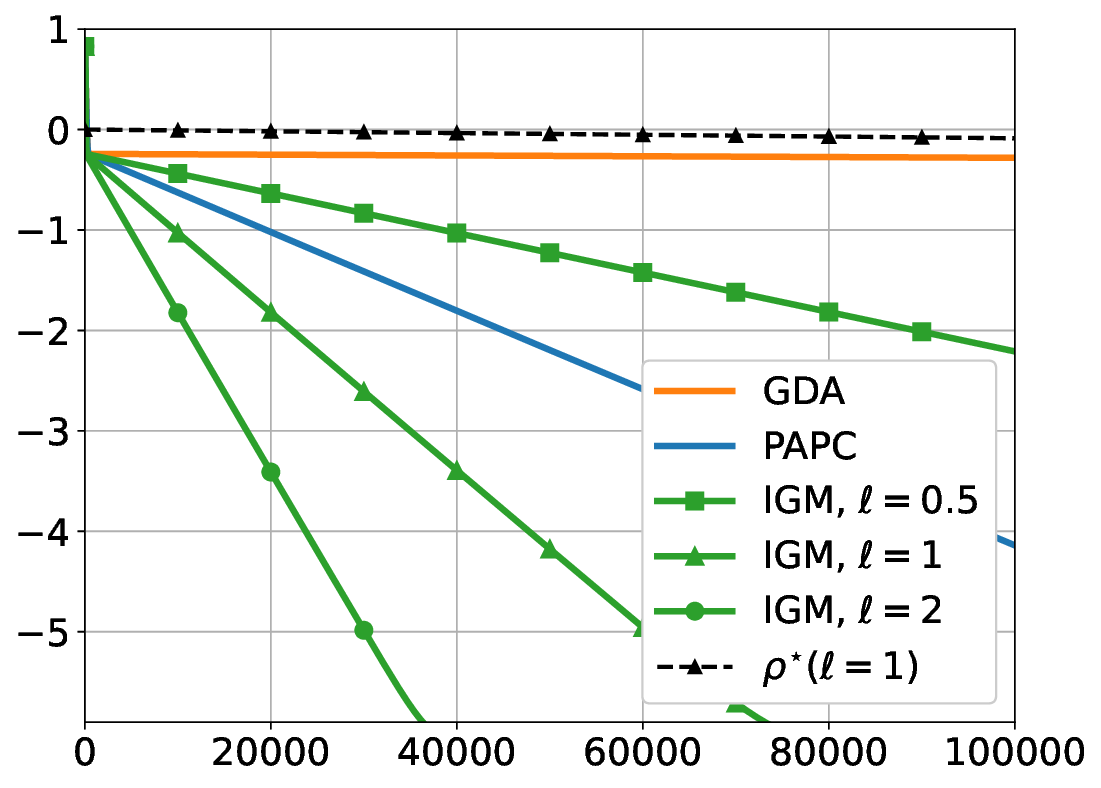}
				\\[-0.0 cm]  {iteration count, $k$}
			\end{tabular}
			&
			\begin{tabular}{c}
				\vspace{.25cm}
				\normalsize{\rotatebox{90}{ $\log_{10} \norm{x^k - \xs}/\norm{\xs}$}}
			\end{tabular}
			&
			\begin{tabular}{c}
				\includegraphics[width=0.4\textwidth]{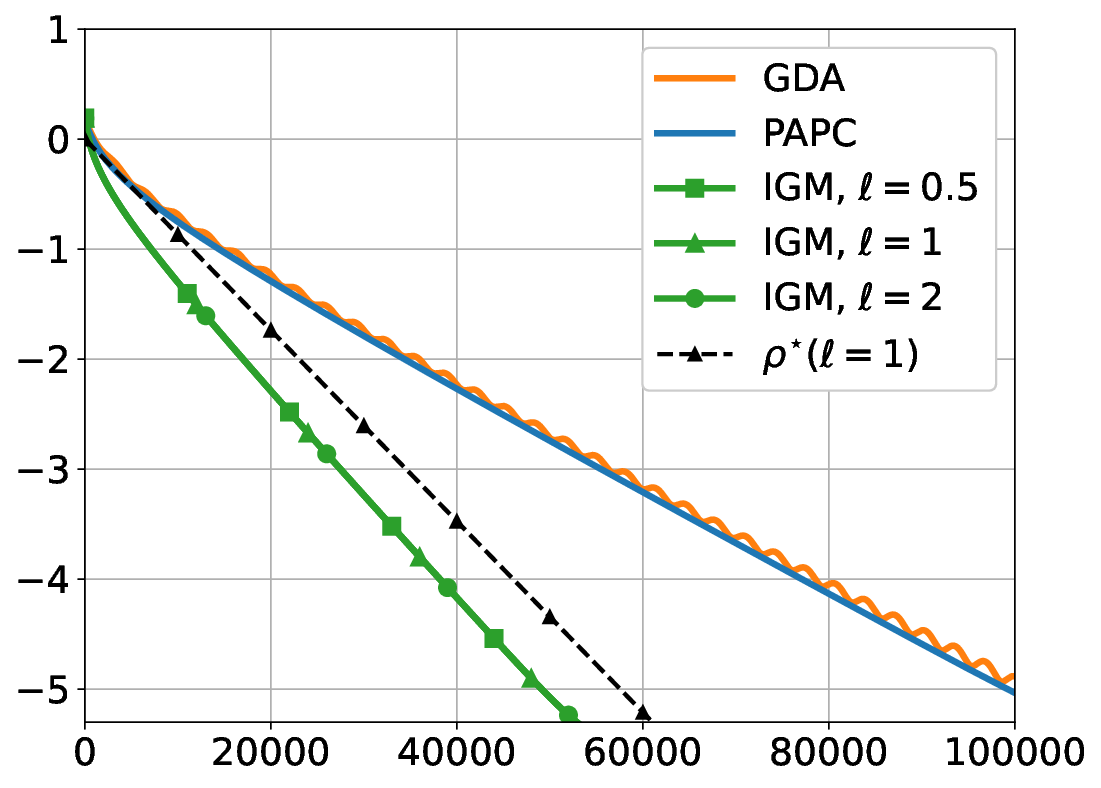}
				\\[-0.0 cm]  {iteration count, $k$}
			\end{tabular}	
		\end{tabular}
		\vspace{0cm}
		\caption{Convergence of the single loop algorithms on Example~1. Left: $(\kappa_f, \kappa_E) = (10^2, 10^6)$; Right: $(\kappa_f, \kappa_E) =(10^4, 10^2)$. The solid lines show convergence of algorithms while the dashed-line $\rho^\star(\ell=1)$ shows the theoretical convergence rate~\eqref{eq.theorem} with $\ell=1$.}
		\label{fig.exp1}
	\end{figure*}

\begin{figure*}[h]
\centering
\begin{tabular}{c@{\hspace{-0.3 cm}}c@{\hspace{0.1 cm}}c@{\hspace{-0.3 cm}}c}
	\begin{tabular}{c}
		\vspace{.25cm}
		\normalsize{\rotatebox{90}{ $\log_{10} \norm{x^k - \xs}/\norm{\xs}$}}
	\end{tabular}
	&
	\begin{tabular}{c}
		\includegraphics[width=0.4\textwidth]{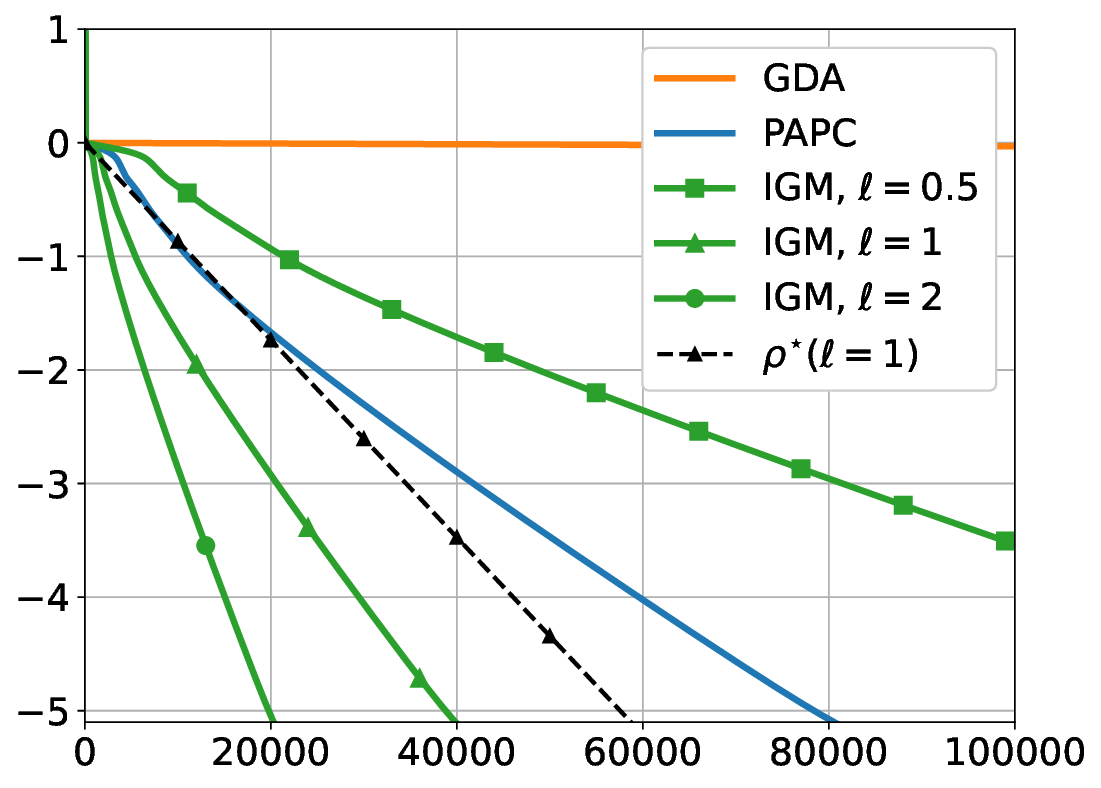}
		\\[-0.0 cm]  {iteration count, $k$}
	\end{tabular}
	&
	\begin{tabular}{c}
		\vspace{.25cm}
		\normalsize{\rotatebox{90}{ $\log_{10} \norm{x^k - \xs}/\norm{\xs}$}}
	\end{tabular}
	&
	\begin{tabular}{c}
		\includegraphics[width=0.4\textwidth]{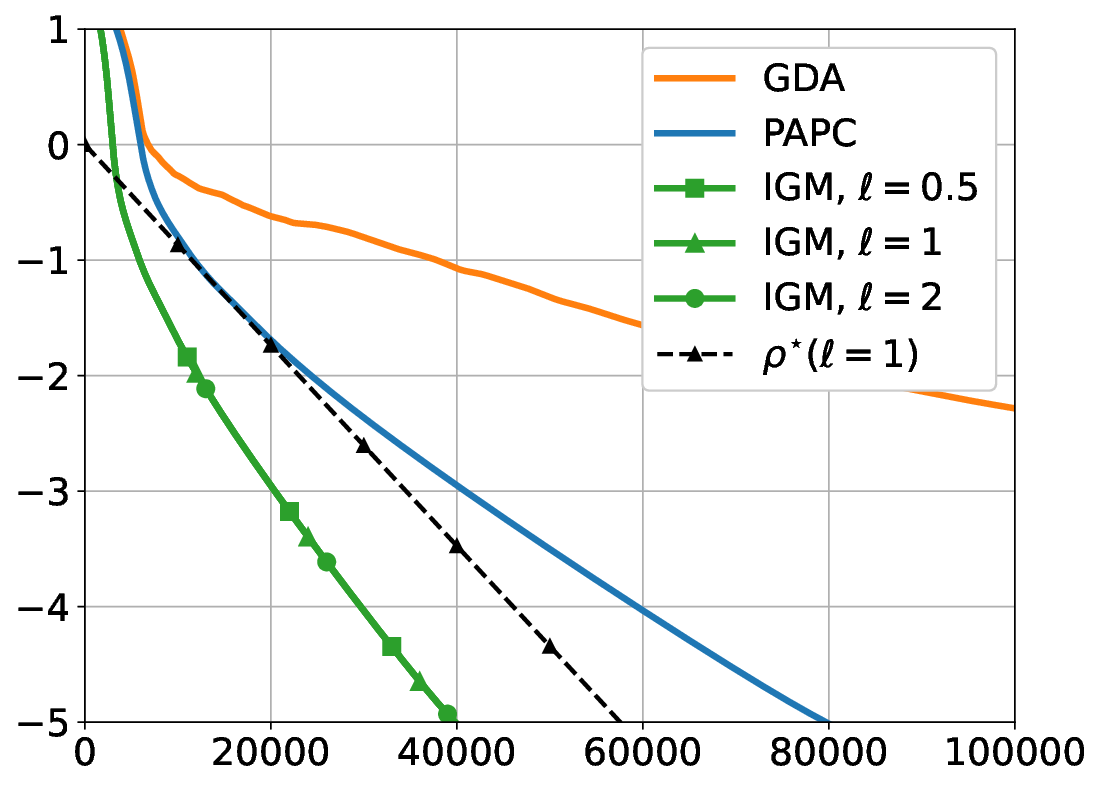}
		\\[-0.0 cm]  {iteration count, $k$}
	\end{tabular}	
\end{tabular}
\vspace{0cm}
\caption{Convergence of the single loop algorithms on Example~2. Left: $(\kappa_f, \kappa_E) = (10^2, 10^4)$; Right: $(\kappa_f, \kappa_E) =(10^5, 10^4)$. The solid lines show convergence of algorithms while the dashed-line $\rho^\star(\ell=1)$ shows the theoretical convergence rate~\eqref{eq.theorem} with $\ell=1$.}
\label{fig.exp2}
\end{figure*}

We demonstrate the validity of our analysis and the effectiveness of the synthesized algorithm through two examples with experimental setups detailed below. For both examples, we compare single-loop gradient-based algorithms that can solve problem~\eqref{eq.main} without requiring any assumptions on constraint matrix $E$ and vector $q$. These include GDA~\cite{algsay20}, PAPC~\cite{salconkovric22} and proposed I-GM (with three different values of $\ell\in\myset{0.5,1,2}$). For each algorithm, we use the parameters that minimize the theoretical convergence rates. Results are presented in Figure~\ref{fig.exp1} and~\ref{fig.exp2}.

While PAPC and I-GM share the same asymptotic complexity, when $\kappa_E\geq \kappa_f$ and \mbox{$\ell=1$}, I-GM requires approximately half the iterations that PAPC needs to reach the same level of accuracy. Consistent with our analysis, increasing $\ell$ reduces the number of I-GM iterations to achieve a given level of accuracy. For instance, increasing $\ell$ from 1 to 2 halves the iteration count, which aligns well with the worst-case convergence rate $\rho$ given in~\eqref{eq.theorem}. Notably, our theoretical lower bound $\rho$ on the convergence rate exhibits a remarkable agreement with the numerical convergence rate, suggesting that the bound is tight.

\paragraph{Example 1:}
We minimize the logistic regression loss,
\beq
f(x) \,=\, \ds \sum_{i=1}^N (-y_ia_i^\top x + \log(1 \,+\, \rme^{a_i^\top x})) + \tfrac{m}{2}\norm{x}_2^2
\eeq
where problem data $A = [a_1 \cdots a_n]^\top \in \reals{N\times n}$ and labels $y_i$'s are generated as described in~\cite{suboycan16}. For the constraint $Ex=q$, we first randomly sample a sparse vector $\xb$ with 50 nonzero entries being set to~1. We then generate constraint matrix $E\in\reals{n\times d}$ with entries sampled from standard normal distribution. After scaling the singular values of $E$ into the desired interval $[\sqrt{\sig_r}, \sqrt{\sig_1}]\cup \myset{0}$, we generate the constraint vector as $q = E\xb$. 

We set the dimensions $(N, n, d)$ to $(10^3, 10^3, 250)$ and the rank of $E$ to $r=200$. We run the algorithms on two instances of this example with different condition numbers. In the first instance, we set $(L,m,\sig_1,\sig_r) = (10, 10^{-1}, 10^4, 10^{-2})$, yielding $\kappa_f=10^2$ and $\kappa_E = 10^6$. In the second instance, we set $(L,m,\sig_1,\sig_r) = (10^2, 10^{-2}, 1, 10^{-2})$, yielding $\kappa_f=10^4$ and $\kappa_E = 10^2$. The results are presented in Figure~\ref{fig.exp1}.

\paragraph{Example 2:}
Inspired from~\cite{salconkovric22}, we conduct a compressed-sensing type experiment. In this experiment, we minimize the cost function,
\beq
f(x) \,=\, \ds \sum_{i=1}^n \bar{f}(x_i),~~~ \bar{f}(x_i) \,=\, \sqrt{x_i^2 \,+\, 1/(L\,-\,m)^2} \,+\, \tfrac{m}{2}\norm{x_i}^2
\eeq
which can be interpreted as an $m$-strongly convex and $L$-smooth approximation of $\ell_1$-norm. The data for equality constraint is generated in the same way as Example~1. 

We set dimensions $(N, n, d)$ to $(10^3, 10^3, 250)$ and the rank of $E$ to $r=200$. Similar to Example~1, we generate two instances. In the first instance, we set $(L,m,\sig_1,\sig_r) = (10, 10^{-1}, 10^2, 10^{-2
})$, yielding $\kappa_f=10^2$ and $\kappa_E = 10^4$. In the latter, we set $(L,m,\sig_1,\sig_r) = (10^4, 10^{-1}, 10^2, 10^{-2})$, yielding $\kappa_f=10^5$ and $\kappa_E = 10^4$. The results are presented in Figure~\ref{fig.exp2}.

\section{Algorithm Representation}
In this section, we demonstrate how formulation~\eqref{eq.repre} encompasses a wide range of algorithms including those with primal-dual updates. As a concrete instance, consider Proximal Alternating Predictor-Corrector (PAPC) algorithm~\cite{salconkovric22}:
\bseq
\begin{align}
x^{k+1/2} & \,=\,  x^k   \,-\, \alpha_1\nabla f(x^k)  \,-\,  \alpha_1 E^\top y^k 
\\
y^{k+1} & \,=\,  y^k  \,+\, \alpha_2(Ex^{k+1/2} \,-\, q)
\label{eq.papc_dual}
\\
x^{k+1} & \,=\,  x^k  \,-\,  \alpha_1\nabla f(x^k) \,-\, \alpha_1 E^\top y^{k+1}.
\end{align}
\eseq		
Replacing the equality constraint with $E^\top Ex = E^\top q$, i.e., multiplying the dual update~\eqref{eq.papc_dual} by $E^\top$, we introduce a new ``dual'' variable $v^k = E^\top y^k$, which reduces the PAPC algorithm to
\bseq\label{eq.prilico2}
\begin{align}
x^{k+1/2} & \,=\,  x^k   \,-\, \alpha_1\nabla f(x^k)  \,-\,  \alpha_1 v^k 
\\
v^{k+1} & \,=\,  v^k  \,+\, \alpha_2E^\top(Ex^{k+1/2} \,-\, q)
\label{eq.papc_dual2}
\\
x^{k+1} & \,=\,  x^k  \,-\,  \alpha_1\nabla f(x^k) \,-\, \alpha_1 v^{k+1}.
\end{align}
\eseq	
Notably, all variables in~\eqref{eq.prilico2} has the same dimension as primal variable $x$. While matching dimensions is not required for our representation~\eqref{eq.repre}, it substantially simplifies the transfer function structure by avoiding block matrix forms. For clarity of exposition, we also assume $q=0$ (this assumption can be avoided by defining an affine operator $W(x) = E^\top(Ex- q)$ as discussed in Section~\eqref{sec.algo}). Substituting the intermediate variable $x^{k+1/2}$ into~\eqref{eq.papc_dual2}, we obtain the following primal-dual updates:
\bseq\label{eq.prilico3}
\begin{align}
v^{k+1} & \,=\,  (I - \alpha_1\alpha_2E^\top E)v^k  \,+\, \alpha_2E^\top E(x^k   \,-\, \alpha_1\nabla f(x^k))
\label{eq.papc_dual3}
\\
x^{k+1} & \,=\,  x^k  \,-\,  \alpha_1\nabla f(x^k) \,-\, \alpha_1 v^{k+1}.
\end{align}
\eseq
In the most general form, any primal-dual update (for solving problem~\eqref{eq.main}) that is linear both in optimization variable and gradients, can be represented as
\bseq\label{eq.repre2}
\begin{align}
z\xh  &\,=\,  \cK_1(z,E)\xh  \,+\, \cK_2(z,E)\gradh   \,+\, \cK_0(z,E)\vh \label{eq.repre2a}
\\
z\vh  &\,=\,  \cL_0(z,E)\vh  \,+\,  \cL_1(z,E)\xh  \,+\, \cL_2(z,E)\gradh.  \label{eq.repre2b}
\end{align}
 \eseq
For instance, the primal version of PAPC algorithm~\eqref{eq.prilico3}, a.k.a., PriLiCo~\cite{salconmisric22}, corresponds to the following transfer function definitions,
\beq\non
\begin{split}
\cK_0(z,E)  &\,=\,  -\alpha_1 z I,
\\
\cK_1(z,E)  &\,=\,  I,
\\
\cK_2(z,E)  &\,=\,  \alpha_1I,
\end{split}~~~~~~
\begin{split}
\cL_0(z,E)  &\,=\,  I  \,-\,  \alpha_1\alpha_2E^\top E
\\
\cL_1(z,E)  &\,=\,  \alpha_2 E^\top E
\\
\cL_2(z,E)  &\,=\,  -\alpha_1\alpha_2 E^\top E.
\end{split}
\eeq
To establish a connection with our representation~\eqref{eq.repre}, we express $\vh$ in terms of  $\xh$ and $\gradh$,
\bseq
\begin{align}
\vh  &\,=\, (zI \,-\, \cL_0)^{-1}(\cL_1\xh  \,+\, \cL_2\gradh)
\\
&\,=\, \widetilde{\cL}_1\xh  \,+\, \widetilde{\cL}_2\gradh\label{eq.vsub}
\end{align}
\eseq
where, omitting the arguments, $\widetilde{\cL}_1\DefinedAs (zI - \cL_0)^{-1}\cL_1$ and $\widetilde{\cL}_2 \DefinedAs (zI - \cL_0)^{-1}\cL_2$. Substituting~\eqref{eq.vsub} into~\eqref{eq.repre2a} yields
\bseq\label{eq.repre3}
\begin{align}
z\xh  &\,=\,  (\cK_1  \,+\,  \cK_0\widetilde{\cL}_1)\xh  \,+\, (\cK_2  \,+\,  \cK_0\widetilde{\cL}_2)\gradh 
\\
&\,=\,  \widetilde{\cK}_1\xh  \,+\,\widetilde{\cK}_2\gradh
\end{align}
\eseq
where $\widetilde{\cK}_1\DefinedAs \cK_1  \,+\,  \cK_0\widetilde{\cL}_1$ and $\widetilde{\cK}_2 \DefinedAs \cK_2  +  \cK_0\widetilde{\cL}_2$. 

The reduction from~\eqref{eq.repre2} to~\eqref{eq.repre3} confirms that algorithms following the primal-dual structure~\eqref{eq.repre2} can be equivalently cast in the primal-only form~\eqref{eq.repre3}, thereby establishing the generality of representation~\eqref{eq.repre}.

\end{document}